\documentclass[12pt]{elsarticle}
\usepackage{amsmath,amssymb,epsfig,amsfonts,epsfig,graphicx}
\usepackage[normalem]{ulem} 
\usepackage{stmaryrd} 
\SetSymbolFont{stmry}{bold}{U}{stmry}{m}{n} 
\usepackage{arydshln} 
\usepackage{epsfig}
\usepackage{tikz-cd}
\usepackage{color}
\usepackage{pict2e}
\usepackage{bbold}

\renewcommand{\phi}{\varphi}
\newcommand{\vertt}[1]{{\left\vert\kern-0.25ex\left\vert\kern-0.25ex\left\vert#1\right\vert\kern-0.25ex\right\vert\kern-0.25ex\right\vert}}

\journal{Mechanics of Materials}

\newtheorem{Remark}{Remark}
\newtheorem{Lemma}{Lemma}
\newtheorem{Definition}{Definition}
\newtheorem{Corollary}{Corollary}
\newtheorem{Theorem}{Theorem}

\renewcommand{\rho}{\varrho}

\newcommand{\ee}{\mathrm{e}}

\newcommand{\dtau}{\,\mathrm{d}\tau}

\newcommand{\D}{\mathrm{D}}
\newcommand{\ddt}{\frac{\mathrm{d}}{\mathrm{d}t}}

\newcommand{\pat}{\partial_t}
\newcommand{\patt}{\partial_{tt}}

\newcommand{\paj}{\partial_j}
\newcommand{\pak}{\partial_k}
\newcommand{\pal}{\partial_l}
\newcommand{\cD}{{\cal D}}
\newcommand{\ve}{\varepsilon}
\newcommand{\ut}{\tilde{u}}
\newcommand{\vet}{\widetilde{\varepsilon}}
\newcommand{\nn}{\nonumber}

\newcommand{\ovu}{\overline{u}}
\newcommand{\ove}{\overline{\varepsilon}}

\newcommand{\R}{\mathbb{R}}
\newcommand{\C}{\mathbb{C}}
\newcommand{\N}{\mathbb N}

\newcommand{\tr}{\mathrm{tr}}

\newcommand{\Div}{\mathrm{Div}}

\newcommand{\sym}{\text{sym}}

\newcommand{\SYM}{\mathrm{Sym}(n)}
\newcommand{\cN}{{\mathcal N}}
\newcommand{\id}{\mathbb{1}}

\DeclareMathOperator\supp{supp}

\begin{document}

\begin{frontmatter}
\title{Global existence and uniqueness of weak solutions\\
for a Willis-type model of elastodynamics}
\author{Thomas Blesgen}
\ead[1]{t.blesgen@th-bingen.de}
\address{Bingen University of Applied Sciences, Berlinstra{\ss}e 109,
D-55411 Bingen, Germany}
\author{Patrizio Neff}
\ead[2]{patrizio.neff@uni-due.de}
\address{Faculty of Mathematics, University of Duisburg-Essen,
Thea-Leymann-Stra{\ss}e 9, D-45127 Essen, Germany}
\cortext[2]{Corresponding author}
\cormark[2]

\begin{abstract}
The existence and uniqueness of 
weak solutions is shown for a system
related to the Willis model of elastodynamics. Both the whole space case and
the case of a bounded smooth domain are studied. To this end the
equations are reformulated as a linear symmetric hyperbolic system of first
order and the existing theory for such systems is applied.
If the initial and boundary data is regular enough, classical solutions
are obtained. The possibility to transform the problem to a linear
symmetric hyperbolic system hinges on a new symmetry condition on the Willis
coupling tensor S, not yet considered in the literature.
This condition demands that $S$ is a totally symmetric third-order tensor.
\end{abstract}

\end{frontmatter}

\section{The Willis model in elastomechanics}
\label{secmodel}

The Willis model, see \cite{Willis81, Willis09, Willis11} and
Milton and Willis \cite{Willis07}, is an extension of classical elastodynamics
with the aim to better reproduce wave propagation in metamaterials.
The latter topic is of high current interest, see, e.g.,
\cite{Rizzi24, Gattin25}.
In fact, metamaterials (or architected materials) show uncommon dispersion
relations that are impossible to predict with classical linear elastodynamics.
In the Willis type models, the balance of linear momentum is modified,
together with a coupling equation. The system reads
\begin{equation}
\label{E1}
\begin{alignedat}{2}
\Div\,\sigma \;=\;& \dot{\mu},\\
\sigma \;=\;& \C_\mathrm{eff}.\sym\,\D u+S_\mathrm{eff}\,\dot{u},\\
\mu \;=\;& S_\mathrm{eff}^T.\sym\,\D u+\rho_\mathrm{eff}\,\dot{u}.
\end{alignedat}
\end{equation}
Moreover, $\C_\mathrm{eff}: \mathrm{Sym}(3)\to\mathrm{Sym}(3)$,
$S_\mathrm{eff}:\R^3\to\mathrm{Sym}(3)\subset\R^{3\times3}$,
$S_\mathrm{eff}^T:\mathrm{Sym}(3)\to\R^3$ have a formal character, but
should be determined by some 'homogenization' procedure. For $S_\mathrm{eff}=0$,
we have just $\dot{\mu}=\rho_\mathrm{eff}\,\ddot{u}$.
The (symmetric) Cauchy stress tensor is denoted by $\sigma$, the
displacement is $u:\Omega\subset\R^3\to\R^3$ and
$\mu:\Omega\subset\R^3\to\R^3$ is the linear momentum density.

In index notation, the system (\ref{E1}) is equivalent to (using summation
convention)
\begin{align}
\label{E1ind1}
\partial_j\sigma_{ij}\;=\;\;& \dot{\mu}_i,\\
\label{E1ind2}
\sigma_{ij} \;=\;\;& (\C_\mathrm{eff})_{ijkl}\;\ve_{kl}
+(S_\mathrm{eff})_{ijk}\;\dot{u}_k,\\
\label{E1ind3}
\mu_i \;=\;\;& (S_\mathrm{eff})_{kli}\;\ve_{kl}+\rho_\mathrm{eff}\,\dot{u}_i.
\end{align}
We show existence and uniqueness of (weak) solutions under the two symmetry
assumptions on the Willis coupling tensor 
\begin{align}
\label{CC1}
S_{ijk} \;=\;& S_{jik}\qquad\mbox{for }1\le i,j,k\le n,\\
\label{CC2}
S_{ijk} \;=\;& S_{jki}\qquad\mbox{for }1\le i,j,k\le n.
\end{align}
While the first symmetry (\ref{CC1}) is naturally associated to
the symmetry of the Cauchy stress $\sigma$, the second condition
(\ref{CC2}) appears to be new and is related to the possibility to
transform the problem to a linear symmetric hyperbolic system.
Combined, the conditions (\ref{CC1}) and (\ref{CC2}) imply that $S$ is a
totally symmetric third-order tensor, having at most 10
independent coefficients in dimension $3$, see \cite{IR25}.
Due to the term $(S_\mathrm{eff}^T)_{lik}=(S_\mathrm{eff})_{kli}$
in (\ref{E1ind3}), the material exhibits reciprocity, see
\cite{Muhl16}, independently of the symmetry assumption~(\ref{CC2}).

\vspace*{2mm}
We would like to thank Prof. G. Rosi (University Paris -- Est Cr{\'e}teil)
for sharing the following observation.
\vspace*{-1mm}
\begin{Remark}
\label{rem1}
Assume that $\rho$, $\C_{ijkl}$ and $S_{ijk}$ are constants independent of
$(x,t)$.
Differentiating (\ref{E1ind2}) w.r.t. $j$ and (\ref{E1ind3}) w.r.t.
$t$, we obtain for $i\in\{1,2,\ldots,n\}$
\begin{align}
\label{R1}
\C_{ijkl}\,\paj\pak u_l+S_{ijk}\,\pat\paj u_k =\;& \paj\sigma_{ij},\\
\label{R2}
S_{jki}\,\pat\paj u_k+\rho\,\pat^2 u_i =\;& \pat\mu_i.
\end{align}
Due to Eq.~(\ref{E1ind1}), the right hand sides of (\ref{R1}), (\ref{R2}) are
equal such that
\begin{equation}
\label{R3}
\rho\,\pat^2 u_i=\C_{ijkl}\,\paj\pak u_l+(S_{ijk}-S_{jki})\,\pat\paj u_k.
\end{equation}
If $S$ satisfies the symmetry (\ref{CC2}), the term in brackets on the
right disappears and Eq.~(\ref{R3}) simplifies to the common linear elasticity
equation
\[ \rho\,\pat^2 u_i=\C_{ijkl}\,\paj\pak u_l \]
showing no Willis coupling. Hence, total symmetry of $S$ excludes Willis
coupling for constant coefficients. However, in general, the coupling tensor
$S$ is space and time dependent in which case the studied Willis problem is
non-trivial in the sense that it departs considerably from classical linear
elasticity. In addition, even for constant $\rho$, $\C$ and $S$, the
Willis coupling does not disappear on $\partial\Omega$ (related to the normal
component $\sigma_{ij}n_j$ of the stress).
Also note that for materials whose underlying microstructure
shows centro-symmetry (inversion symmetry), the third-order Willis coupling
tensor $S$ must vanish. This is not at odds with assuming that $S$ is totally
symmetric as both conditions are different. 
\end{Remark}

\subsection{Background for the Willis approach}
Unusual dynamic properties of certain classes of composite materials
(architected materials, metamaterials) necessitate to introduce new systems
of equations, extending classical linear elasticity. One direction is to
explore generalized continua, see. e.g. \cite{Mad16,Rizzi24}, with additional
kinematic descriptor fields. Another direction is to change the structure of
the equations. In the latter case, the Willis model aims to provide
effective constitutive equations for ensemble averages of quantities of
interest, i.e. Cauchy stress $\sigma$ versus displacement $u$ or velocity
$v=\dot{u}$.
Hence, the usefulness of the Willis approach depends crucially on the
assumption that an ensemble average is a reasonable descriptor of the given
(periodic) microstructure of the medium.

Let the fully resolved microstructure of the medium obey classical
linear elasticity, i.e.
\begin{equation}
\label{L1}
\Div\,\sigma \;=\; \dot{\mu} \;=\; \ddt(\rho\,\dot{u})
\end{equation}
together with the constitutive law
\begin{equation}
\label{L2}
\sigma \;=\; \C(x).\ve,\qquad \ve \;=\; \sym\, \D u.
\end{equation}
Ensemble averaging (\ref{L1}) we obtain
\[ \Div\,\langle\sigma\rangle \;=\; \langle\dot{\mu}\rangle \]
for the averaged quantities $\langle\sigma\rangle$ and
$\langle\dot{\mu}\rangle$. The infinitesimal strain tensor is
likewise ensemble averaged as
\begin{equation}
\label{L3}
\langle\ve\rangle \;=\; \sym\, \D\langle u\rangle.
\end{equation}
However, ensemble averaging the constitute law (\ref{L2})$_1$ is not directly
achievable since
\begin{equation}
\label{L4}
\langle\sigma\rangle \;=\; \langle\C(x).\ve\rangle \;\not=\;
\langle\C(x)\rangle.\langle\ve\rangle
\end{equation}
and $\langle\mu\rangle \not= \langle\rho\rangle\cdot\langle\dot{u}\rangle$
because in general the product of averages differs from the average of
the product. Therefore, the Willis equations provide simple constitutive
closure relations in the form
\begin{equation}
\label{L5}
\begin{aligned}
\langle\sigma\rangle \;:=\;& \C_\mathrm{eff}.\langle\ve\rangle
+S_\mathrm{eff}.\langle\dot{u}\rangle,\\
\langle\mu\rangle \;:=\;& S_\mathrm{eff}^T.\langle\ve\rangle
+\rho_\mathrm{eff}\,\langle\dot{u}\rangle,
\end{aligned}
\end{equation}
where $\C_\mathrm{eff}$, $S_\mathrm{eff}$, $\rho_\mathrm{eff}$ must be
determined in an additional step. In the following, we skip the
$\langle\cdot\rangle$-notation.

\subsection{Invariance considerations}
If we adhere to the idea that the Willis system (\ref{E1}) should be the result
of some homogenization based on classical linear elastodynamics, then it is
natural to require that solutions of (\ref{E1}) should satisfy
the same invariance conditions as are satisfied by classic linear elastodynamics
\begin{equation}
\label{E2}
\Div\,\sigma\;=\; \rho\, \patt u,\qquad \sigma \;=\; \C.\sym\,\D u.
\end{equation}
It is easy to see that (\ref{E2}) is {\it infinitesimal Galilean-invariant},
i.e. if $u$ is a solution, so is
\begin{equation}
\label{E3}
u(x)\mapsto u(x)+\overline{A}x+\overline{b}(t),\qquad
\overline{b}''(t)=0,
\end{equation}
where $\overline{A}\in\mathfrak{so}(3)$ and $t\mapsto\overline{b}(t)\in\R^3$.
A direct check reveals that the system~(\ref{E1}) likewise
admits the invariance (\ref{E3}).
\footnote{It is clear that an ensemble average (statistical average)
transforms likewise.}

However, the system~(\ref{E2}) also admits a lesser-known further
invariance condition, the so-called {\it extended infinitesimal Galilean
invariance}
\begin{equation}
\label{E4}
u(x)\mapsto u(x)+\overline{A}(t)x+\overline{b}(t),\qquad
\overline{b}''(t)=0,\quad \overline{A}''(t)=0.
\end{equation}
This invariance condition has no immediate counterpart in classical nonlinear
elasto-dynamics but appears as possibility due to the loss of information
inherent in the linearization process of which (\ref{E2}) is the
result. Be that as it may, it must be observed
that the linear Willis system is not invariant w.r.t. (\ref{E4}) if
$S_\mathrm{eff}\not=0$.
Thus, whatever process of homogenization is applied, the simplified system
(\ref{E1}) cannot entirely capture all effects that are possible in
a fully dynamic calculation of a completely resolved microstructure.
Nevertheless, we find it worthwhile to look at the mathematical structure
presented by the system~(\ref{E1}).
To the best of our knowledge, no local or global existence proof has yet been
given.
Due to the linearity, however, this should be possible (but see \cite{Lewy57})
and indeed, based on the general theory of linear symmetric hyperbolic
systems of first order, the abstract Willis system can be cast into a format
that permits an existence result.

\section{Prerequisites and assumptions}
\label{secpre}
Let $I:=\{1,2,\ldots,n\}$ and $\Omega\subset\R^n$ be a domain, $0<T\le\infty$
a fixed time,
$\Omega_T:=\Omega\times(0,T)$ and $\cD:=\Omega\times[0,T]$.
If $\Omega$ is bounded we write $\Sigma_T:=\partial\Omega\times(0,T)$.

Throughout, we shall employ the following notations.
We write $\partial_k$ shortly for $\frac{\partial}{\partial x_k}$ and
$\|M\|:=\tr(M^TM)$ is the Frobenius norm, where $\tr(M):=\sum_{k\in I}M_{kk}$
is the trace and $M^T$ the transpose of $M$ for $M\in\R^{n\times n}$.
We write $\langle v,w\rangle:=\sum_{k\in I}v_kw_k$ for the Euclidean inner
product of two vectors $v,w\in\R^n$ and $\C.\ve:=\C_{ijkl}\,\ve_{kl}$ for the
application of a fourth-order tensor $\C$ to a second-order tensor $\ve$.
Let $s\in\N$ be a fixed integer. By $W^{s,p}(\Omega)$ we denote the Sobolev
space of $s$-times weakly differentiable functions in $L^p(\Omega)$ and
$H^s(\Omega)\equiv W^{s,2}(\Omega)$ is a Hilbert space.
By $C_b^m(X)$ we denote the space of $m$-times bounded differentiable functions
of a Banach space $X$ to $\R$. By $\SYM$ we denote the set of symmetric real
$n\times n$ matrices.

For $(x,t)\in\Omega_T$, let
$u=u(x,t)=(u_i(x,t))_{i\in I}$ denote the displacement,
$\mu=\mu(x,t)=(\mu_i(x,t))_{i\in I}$ the momentum density vector,
$\sigma=\sigma(x,t)=(\sigma_{ij})_{i,j\in I}$ the symmetric
Cauchy stress tensor; $S=S(x,t)=(S_{ijk}(x,t))_{i,j,k\in I}$ is the Willis
coupling tensor, $\pat u$ is the particle velocity.

\vspace*{2mm}
One model assumption is that the strain be small.
By $\ve=\ve(u):=\sym\,\D u$ we denote the linearized strain tensor, i.e.
\begin{equation}
\label{vedef}
\ve_{kl}:=\frac12\Big(\pak u_l+\pal u_k\Big),\qquad k,l\in I.
\end{equation}
For the existence proofs below we make the following assumptions on
$\rho$, $\C$, $S$, $\ovu$ and $\mu_0$.

\vspace*{2mm}
\noindent{\bf(A0)} The elasticity tensor $\C$ is a fourth-order tensor with
$(\C=\C_{ijkl}(x,t))_{i,j,k,l\in I}$ possibly depending on $(x,t)$ to account
for complicated material behaviour. We assume
\begin{equation}
\label{CCb}
\C_{ijkl}(x,t)\in C_b^\infty(\Omega\times[0,T]).
\end{equation}
The tensor $\C$ satisfies the {\it major} and {\it minor symmetry relations}
\begin{equation}
\label{Csym}
\C_{ijkl}(x,t)=\C_{jikl}(x,t)=\C_{ijlk}(x,t)=\C_{klij}(x,t)
\qquad\mbox{for all }i,j,k,l\in I,\, (x,t)\in\cD.
\end{equation}
We assume that $\C$ is uniformly positive definite.
This means there exists a constant $c_1>0$ such that for all $(x,t)\in\cD$
\begin{equation}
\label{PD}
\big\langle\C(x,t).\ve,\ve\big\rangle \;\ge\; c_1\|\ve\|^2
\qquad\mbox{for all }\ve\in\SYM.
\end{equation}

\noindent{\bf(A1)} The mass density $\rho=\rho(x,t)$ of the material is
given and satisfies
\begin{equation}
\label{rhoCb}
\rho,\,\pat\rho\in C_b^\infty(\Omega\times[0,T]).
\end{equation}
There exists a constant $m_0>0$ such that
\begin{equation}
\label{rhopos}
\rho(x,t)\ge m_0\qquad\mbox{for all }(x,t)\in\Omega\times[0,T].
\end{equation}
\vspace*{2mm}
\noindent{\bf(A2)} The third-order Willis coupling tensor
$S=S_{ikl}(x,t)$ satisfies
\begin{equation}
\label{SCb}
S_{ikl},\,\pat S_{ikl}\in C_b^\infty(\Omega\times[0,T])\qquad
\mbox{for all }i,k,l\in I.
\end{equation}
The tensor $S$ satisfies the symmetry relations
\begin{align}
\label{Ssym}
S_{ijk}(x,t)=\;& S_{jik}(x,t)\qquad\mbox{for all }i,j,k\in I,
(x,t)\in\Omega\times[0,T],\\
\label{Ssym2}
S_{ijk}(x,t)=\;& S_{jki}(x,t)\qquad\mbox{for all }i,j,k\in I,
(x,t)\in\Omega\times[0,T].
\end{align}
\vspace*{2mm}
\noindent{\bf(A3)} The initial data $u_0$ and $\mu_0$ satisfy for an integer
$s\ge1$
\begin{align}
\label{u0reg}
u_0\in\;& H^{s+1}(\Omega;\,\R^n),\\
\label{mu0reg}
\mu_0\in\;& H^s(\Omega;\,\R^n).
\end{align}
\vspace*{2mm}
\noindent{\bf(A4)} The boundary function $\ovu$ can be extended to a function
on $\overline{\Omega}\times[0,T]$ which satisfies
\begin{align}
\label{A4C3a}
\pat^r\,\ovu(\cdot,0)\in H^{s+1-r}(\Omega;\,\R^n)
\qquad\mbox{for }\/0\le r\le s+1,\\
\label{A4C3b}
\pat^r\,\ovu\in L^2(0,T;\,H^{s+2-r}(\Omega;\,\R^n))
\qquad\mbox{for }\/0\le r\le s+2.
\end{align}

\vspace*{5mm}
\noindent
We write $\rho_0(x):=\rho(x,0)$ for the (given) density at time $t=0$.
The condition~(\ref{Ssym}) ensures the symmetry of the Cauchy stress
tensor $\sigma$.
The boundary data $\ovu$ in (A4) is introduced below in (\ref{split}).
We assume the compatibility of initial and boundary data, i.e.
\begin{equation}
\label{compa}
\ovu(\cdot,0)=u_0\qquad\mbox{in }\Omega.
\end{equation}

\vspace*{2mm}
Due to (\ref{Csym}), for $n=3$, only $36$ of the $81$ entries of $\C$ are
independent. As is well known, see \cite[pp. 268--269]{Sommerfeld}, due to
conservation of energy, this reduces further and at most $21$ entries may be
independent. The material symmetry relation allows to reduce this number even
more, see, e.g., \cite{MC90,Van18}.

Using the symmetries (\ref{Csym}), we recover the identity
\begin{align}
(\C(x,t).\ve)_{ij} \;=\;\; & \C_{ijkl}(x,t)\,\frac12\big(\pak u_l+\pal u_k\big)
\;=\; \frac12\C_{ijkl}(x,t)\,\pak u_l+\frac12\C_{ijkl}(x,t)\,\pal u_k\nn\\
\label{epsu}
=\;\;& \frac12\big(\C_{ijkl}(x,t)+\C_{ijlk}(x,t)\big)\pak u_l\nn\\
=\;\;& \C_{ijkl}(x,t)\,\pak u_l \;=\; (\C(x,t).\D u)_{ij},\qquad i,j\in I,\;
(x,t)\in\cD.
\end{align}

\vspace*{1mm}
Subsequently we analyze the following system of equations related to the
{\em Willis model}.

\noindent Find the solution vector $(\mu,\sigma,u)$ with
$\mu\in L^2(0;T; L^2(\Omega;\,\R^n))$,
$u\in L^2(0,T;\,H^1(\Omega;\,\R^n))$,
$\sigma\in L^2(0;T;\, H^1(\Omega;\,\R^{n\times n}))$ solving in $\Omega_T$
\begin{align}
\label{S1o}
\pat\mu_i \;=\;\;& \paj\sigma_{ij},\hspace*{65pt} i\in I,\\
\label{S2o}
S_{ijk}\,\pat u_k \;=\;\;& \sigma_{ij}-\C_{ijkl}\,\ve_{kl},\qquad i,j\in I,\\
\label{S3o}
\rho\,\pat u_i \;=\;\;& \mu_i-S_{kli}\,\ve_{kl},\hspace*{32pt}i\in I
\end{align}
subject to the initial and boundary conditions
\begin{align}
\label{BC1}
u(\cdot,0) \;=\;\; & u_0, \qquad\mbox{in }\Omega,\\
\label{BC2}
\mu(\cdot,0) \;=\;\; & \mu_0, \qquad\mbox{in }\Omega,\\
\label{BC3}
u \;=\;\;& \ovu, \hspace*{29pt}\mbox{on }\partial\Omega\times[0,T]
\end{align}
for given initial values $u_0\in H^{1,2}(\Omega)$ and $\mu_0\in L^2(\Omega)$.
In (\ref{epsu}), (\ref{S1o})--(\ref{S3o}) and below, we utilize the summation
convention and implicitly sum over repeated indices in $I$ unless stated
otherwise. In the original Willis model \cite{Willis85}, formulated in $n=3$
space dimensions, $S$ is defined by a convolution. In this article, we do not
assume any specific form of $S$, but consider generic tensors $S$ depending on
$(x,t)$.

\begin{Remark}
\label{rem2}
For $S={\bf0}$, (\ref{S1o})-(\ref{S3o}) constitute the classical
equations of motion for the propagation of waves in solids and from
(\ref{S2o}) we recover {\em Hooke's law}
\begin{equation}
\label{Hooke} \sigma=\C.\ve.
\end{equation}
\end{Remark}

\section{Reformulation of the problem}
\label{secreform}
We assume that the boundary data $\ovu$ in (\ref{BC3}) can be extended
to a function $\ovu\in H^{s+2}(\overline{\Omega}\times[0,T];\,\R^n)$.
We split the deformation vector $u$ by writing
\begin{equation}
\label{split}
u(x,t)\,:=\,\ut(x,t)+\ovu(x,t)
\end{equation}
such that $\ut=0$ on $\partial\Omega\times[0,T]$.
Due to compatibility of initial and boundary data, $\ovu(\cdot,0)-u_0$
has zero trace on $\partial\Omega$.
Analogous to the definition of $\ve=\ve(u)$ in (\ref{vedef}), we set
\[ \vet\,:=\,\ve(\ut)\in\SYM,\qquad \ove\,:=\,\ve(\ovu)\in\SYM. \]
With these notations, (\ref{S1o})--(\ref{S3o}) rewrites as the following
system in $\Omega_T$
\begin{align}
\label{S1}
\pat\mu_i \;=\;\;& \paj\sigma_{ij},\hspace*{103pt} i\in I,\\
\label{S2}
S_{ijk}\,\pat(\ut_k+\ovu_k) \;=\;\;& \sigma_{ij}-\C_{ijkl}\,
(\vet_{kl}+\ove_{kl}), \qquad i,j\in I,\\
\label{S3}
\rho\,\pat(\ut_i+\ovu_i) \;=\;\;& \mu_i-S_{kli}\,
(\vet_{kl}+\ove_{kl}),\hspace*{32pt}i\in I
\end{align}
subject to the initial and boundary conditions (\ref{BC1})--(\ref{BC3}).

We reformulate (\ref{S1})--(\ref{S3}) as a linear hyperbolic system.
Using (\ref{rhopos}), Eqn.~(\ref{S3}) becomes
\begin{equation}
\label{s1}
\pat(\ut_m+\ovu_m) \;=\; \frac{1}{\rho}\mu_m-\frac{1}{\rho}S_{klm}\,
\big(\vet_{kl}+\ove_{kl}\big),\qquad m\in I.
\end{equation}
\vspace*{-2.5mm}
Hence
\begin{align}
\sigma_{ij} \stackrel{(\ref{S2})}{\;=\;}& \C_{ijkl}\,\big(\vet_{kl}+\ove_{kl}
\big)+S_{ijm}\,\pat\big(\ut_m+\ovu_m\big)\nn\\
\label{Hdef}
\stackrel{(\ref{s1})}{=\;}&
\underbrace{\Big(\C_{ijkl}-\frac{1}{\rho}S_{ijm}S_{klm}
\Big)}_{=:\; H_{ijkl}}
\big(\vet_{kl}+\ove_{kl}\big)+\frac{1}{\rho}S_{ijk}\,\mu_k,\qquad i,j\in I.
\end{align}
Plugging (\ref{Hdef}) into (\ref{S1}), we obtain
\begin{equation}
\label{step1}
\pat\mu_i \;=\; \paj\Big(H_{ijkl}\,\big(\vet_{kl}+\ove_{kl}\big)
+\frac{1}{\rho}S_{ijk}\,\mu_k\Big),\qquad i\in I.
\end{equation}
By (\ref{S3}), $\mu_i=\rho\,\pat\big(\ut_i\!+\!\ovu_i\big)
+S_{kli}\,\big(\vet_{kl}\!+\!\ove_{kl}\big)$. 
Using this relationship on the left of (\ref{step1}) yields
\begin{align*}
\dotuline{\pat\big(\rho\,\pat(\ut_i+\ovu_i)+S_{kli}\,
(\vet_{kl}+\ove_{kl})\big)}
=\;& \dashuline{\paj H_{ijkl}\,\big(\vet_{kl}+\ove_{kl}\big)}
+\uuline{H_{ijkl}\,\paj\big(\vet_{kl}+\ove_{kl}\big)}\\
& +\uwave{\paj(\rho^{-1})S_{ijk}\,\mu_k}
+\frac{1}{\rho}\paj S_{ijk}\,\mu_k+\uline{\frac{1}{\rho}S_{ijk}\,\paj\mu_k}.
\end{align*}
Plugging in
$\mu_k=\rho\,\pat(\ut_k+\ovu_k)+S_{mnk}\,
\big(\vet_{mn}+\ove_{mn}\big)$ on the right
together with the definition (\ref{Hdef}) of $H$, we obtain
\begin{align}
\dotuline{\rho\,\patt\ut_i}\;=\;&
\dashuline{-\paj(\rho^{-1})S_{ijm}S_{klm}\,\big(\vet_{kl}\!+\!\ove_{kl}\big)
-\frac{1}{\rho}\paj S_{ijm}S_{klm}\,\big(\vet_{kl}\!+\!\ove_{kl}\big)
-\frac{1}{\rho}S_{ijm}\paj S_{klm}\,\big(\vet_{kl}\!+\!\ove_{kl}\big)}\nn\\
&+\dashuline{\paj\C_{ijkl}\,\big(\vet_{kl}+\ove_{kl}\big)}
+\uuline{\Big(\C_{ijkl}-\frac{1}{\rho}S_{ijm}S_{klm}\Big)
\paj\big(\vet_{kl}+\ove_{kl}\big)}\nn\\
& +\uwave{\rho\,\paj(\rho^{-1})\,S_{ijk}\,\pat\big(\ut_k+\ovu_k\big)
+\paj(\rho^{-1})S_{ijk}S_{mnk}\,\big(\vet_{mn}+\ove_{mn}\big)}\nn\\
\label{step2}
& +\paj S_{ijk}\,\pat\big(\ut_k+\ovu_k\big)+\frac{1}{\rho}
\paj S_{ijk}S_{mnk}\,\big(\vet_{mn}+\ove_{mn}\big)\\
& +\uline{\frac{1}{\rho}S_{ijk}\,\paj\rho\,\pat\big(\ut_k+\ovu_k\big)
+S_{ijk}\,\paj\big(\pat\ut_k+\pat\ovu_k\big)
+\frac{1}{\rho}S_{ijk}\paj S_{mnk}\,\big(\vet_{mn}+\ove_{mn}\big)}\nn\\
& +\uline{\frac{1}{\rho}S_{ijk}S_{mnk}\,\paj\big(\vet_{mn}+\ove_{mn}\big)}\nn\\
& -\dotuline{\pat S_{kli}\,\big(\vet_{kl}+\ove_{kl}\big)
-S_{kli}\,\pat\big(\vet_{kl}+\ove_{kl}\big)
-\pat\rho\,\pat\big(\ut_i+\ovu_i\big)-\rho\,\patt\ovu_i},\qquad i\in I.\nn
\end{align}
Most terms in (\ref{step2}) cancel out. After simplifications, we are left with
\begin{align}
\rho\,\patt\ut_i \;=\;& \C_{ijkl}\,\paj\big(\vet_{kl}+\ove_{kl}\big)
+\paj S_{ijk}\,\pat\big(\ut_k+\ovu_k\big)
+(S_{ijk}-S_{jki})\,\pat\paj\big(\ut_k+\ovu_k\big)\nn\\
\label{step3}
& +\paj\C_{ijkl}\,\big(\vet_{kl}+\ove_{kl}\big)
-\pat S_{kli}\,\big(\vet_{kl}+\ove_{kl}\big)
-\pat\rho\,\big(\pat\ut_i+\pat\ovu_i\big)-\rho\,\patt\ovu_i.
\end{align}

\noindent
With (\ref{epsu}), we have
\[ \C_{ijkl}\,\paj\big(\vet_{kl}+\ove_{kl}\big) \;=\; \C_{ijkl}\,\paj\frac12
\big(\pak\ut_l+\pal\ut_k+\pak\ovu_l+\pal\ovu_k\big) \;=\;
\C_{ijkl}\,\paj\pak\big(\ut_l+\ovu_l\big). \]
Similarly, as a consequence of (\ref{Ssym}), the third term on the
right of (\ref{step3}) disappears, and
\begin{align}
\label{dis1}
\paj\C_{ijkl}\,\big(\vet_{kl}+\ove_{kl}\big) \;=\;\;& \paj\C_{ijkl}\,
\pak\big(\ut_l+\ovu_l\big),\\
\label{dis2}
\pat S_{kli}\, \big(\vet_{kl}+\ove_{kl}\big) \;=\;\;& \pat
S_{kli}\, \pak\big(\ut_l+\ovu_l\big).
\end{align}
Eventually, after introducing $e:\cD\to\R^n$ by
\begin{align}
\rho\, e_i \;:=\;& -\C_{ijkl}\,\paj\pak\ovu_l-\paj S_{ijk}\,\pat\ovu_k
+(\pat S_{kli}-\paj\C_{ijkl})\,\pak\ovu_l\nn\\
\label{edef}
& +\pat\rho\,\pat\ovu_i+\rho\,\patt\ovu_i+(S_{jki}-S_{ijk})\pat\paj
\ovu_k,\qquad i\in I
\end{align}
we end up with
\begin{align}
\rho\,\patt\ut_i \;=\;& \C_{ijkl}\,\paj\pak\ut_l+\big(\paj\C_{ijkl}
-\pat S_{kli}\big)\,\pak\ut_l+\paj S_{ijk}\,\pat\ut_k\nn\\
\label{mod}
& -\pat\rho\,\pat\ut_i+(S_{ijk}-S_{jki})\paj\pat\ut_k
-\rho e_i\hspace*{70pt}\mbox{in }\Omega_T,\; i\in I.
\end{align}
The system~(\ref{mod}) is solved subject to the initial and boundary conditions
\begin{align}
\label{mod1}
\rho_0\,\pat\ut(x,0) \;=\;\;& g(x),\qquad x\in\Omega,\\
\label{mod2}
\ut(x,0) \;=\;\;& h(x),\qquad x\in\Omega,\\
\label{mod3}
\ut(x,t) \;=\;\;& 0,\hspace*{40pt} x\in\partial\Omega,\,t\in[0,T].
\end{align}
Therein, the initial data $h$ is specified from (\ref{split}) and (\ref{BC1}),
$g$ is specified from (\ref{S3}) at $t=0$,
\begin{align}
\label{BC1n}
h(x) \;=\;& u_0(x)-\ovu(x,0),\qquad x\in\Omega,\\
\label{BC2n}
g(x) \;=\;& \Big(\mu_{0i}(x)-S_{kli}(x,0)\,\pak\big(u_{0l}(x)+\ovu_l(x,0)\big)
-\rho_0(x)\,\pat\ovu_i(x,0)\Big)_{1\le i\le n},\quad x\in\Omega.
\end{align}
With (\ref{u0reg}), it holds $h\in H_0^{s+1}(\Omega)$ due to
the compatibility of initial and boundary data.

The equations~(\ref{mod}) constitute a linear hyperbolic system and represent
the most general form the resulting equations may have under the assumption
(\ref{Ssym}).

We recall that a {\em first order linear symmetric hyperbolic system}
is of the form
\begin{equation}
\label{linsym}
Lv:=A_0(x,t)\,\pat v+\sum_{k=1}^n A_k(x,t)\,\pak v+B(x,t)\,v=w(x,t),
\end{equation}
where $v:\Omega\times\R_{\ge0}\to\R^m$, $A_0(x,t)$,
$A_1(x,t),\ldots,A_n(x,t)\in\R^{m\times m}$ are symmetric matrices for all
$(x,t)\in\overline{\Omega}\times[0,T]$, $A_0(x,t)$ is positive definite,
$B(x,t)\in\R^{m\times m}$, and $w(x,t)\in\R^m$ is a given right hand side.
As a consequence of
the symmetry of $A_0$, it is diagonalizable and all eigenvalues are real.

As a preparation of the following lemma, we introduce the symmetric $3\times3$
matrices
{\small
\begin{align}
\C_1^1:=\;& \left(\begin{smallmatrix}
\C_{1111}\; & \C_{1112}\; & \C_{1131}\\
\C_{2111}\; & \C_{2112}\; & \C_{2131}\\
\C_{3111}\; & \C_{3112}\; & \C_{3131}
\end{smallmatrix}\right),\quad
\C_2^1:=\left(\begin{smallmatrix}
\C_{1211}\; & \C_{1221}\; & \C_{1231}\\
\C_{2211}\; & \C_{2221}\; & \C_{2231}\\
\C_{3211}\; & \C_{3221}\; & \C_{3231}
\end{smallmatrix}\right),\quad
\C_3^1:=\left(\begin{smallmatrix}
\C_{1311}\; & \C_{1321}\; & \C_{1331}\\
\C_{2311}\; & \C_{2321}\; & \C_{2331}\\
\C_{3311}\; & \C_{3321}\; & \C_{3331}
\end{smallmatrix}\right),\nn\\
\label{Cijdef}
\C_1^2:=\;& \left(\begin{smallmatrix}
\C_{1112}\; & \C_{1122}\; & \C_{1132}\\
\C_{2112}\; & \C_{2122}\; & \C_{2132}\\
\C_{3112}\; & \C_{3122}\; & \C_{3132}\\
\end{smallmatrix}\right),\quad
\C_2^2:=\left(\begin{smallmatrix}
\C_{1212}\; & \C_{1222}\; & \C_{1232}\\
\C_{2212}\; & \C_{2222}\; & \C_{2232}\\
\C_{3212}\; & \C_{3222}\; & \C_{3232}
\end{smallmatrix}\right),\quad
\C_3^2:=\left(\begin{smallmatrix}
\C_{1312}\; & \C_{1322}\; & \C_{1332}\\
\C_{2312}\; & \C_{2322}\; & \C_{2332}\\
\C_{3312}\; & \C_{3322}\; & \C_{3332}
\end{smallmatrix}\right),\\
\C_1^3 :=\;& \left(\begin{smallmatrix}
\C_{1113}\; & \C_{1123}\; & \C_{1233}\\
\C_{2113}\; & \C_{2123}\; & \C_{2233}\\
\C_{3113}\; & \C_{3123}\; & \C_{3233}
\end{smallmatrix}\right),\quad
\C_2^3:=\left(\begin{smallmatrix}
\C_{1213}\; & \C_{1223}\; & \C_{1233}\\
\C_{2213}\; & \C_{2223}\; & \C_{2233}\\
\C_{3213}\; & \C_{3223}\; & \C_{3233}
\end{smallmatrix}\right),\quad
\C_3^3:=\left(\begin{smallmatrix}
\C_{1313}\; & \C_{1323}\; & \C_{1333}\\
\C_{2313}\; & \C_{2323}\; & \C_{2333}\\
\C_{3313}\; & \C_{3323}\; & \C_{3333}
\end{smallmatrix}\right).\nn
\end{align}
}
It holds $C_i^j=C_i^j(x,t)$ with the symmetry $\C_i^j(x,t)=\C_j^i(x,t)$
for $i,j\in\{1,2,3\}$.

For the second term on the right of (\ref{mod}),
we introduce the short-hand notation
\begin{equation}
\label{Ddef}
D_{kl}^i(x,t):=\pat S_{kli}(x,t)-\paj\C_{ijkl}(x,t),\qquad i,k,l\in\{1,2,3\}.
\end{equation}
If $\Omega$ is bounded and $\nu=(\nu_1,\nu_2,\nu_3)\in\R^3$ is the
unit outer normal vector at a point in $\partial\Omega$, let
\begin{equation}
\label{Cnudef}
\C_\nu:=\left(\begin{array}{c;{2pt/2pt}c;{2pt/2pt}c}
\nu_1\C_1^1\; & \nu_1\C_2^1\; & \nu_1\C_3^1\\ \hdashline[2pt/2pt]
\nu_2\C_1^2\; & \nu_2\C_2^2\; & \nu_2\C_3^2\\ \hdashline[2pt/2pt]
\nu_3\C_1^3\; & \nu_3\C_2^3\; & \nu_3\C_3^3\end{array}\right)
\!\in\R^{9\times9}.
\end{equation}
\begin{Lemma}
\label{lem1}
Assume that the Willis coupling tensor satisfies the further symmetry relation
\begin{equation}
\label{Ssym2b}
S_{ijk}(x,t)=S_{jki}(x,t)\qquad\mbox{for all }i,j,k\in I,\,
(x,t)\in\Omega\times[0,T].
\end{equation}
Then the equations (\ref{mod})--(\ref{mod3}) constitute a linear symmetric
hyperbolic system of first order, i.e. they can be written as the mixed
initial boundary value problem
\begin{align}
\label{syst1}
Lv &= w \hspace*{42pt} \mbox{in }\Omega\times[0,T],\\
\label{syst2}
M(x)v &=0 \hspace*{45pt} \mbox{in }\Gamma\times[0,T],\\
\label{syst3}
v(\cdot,0) &= v_0 \hspace*{40pt}\mbox{in }\Omega,
\end{align}
where $L$ is given by (\ref{linsym}), $M(x)\in\R^{m\times m}$ for $x\in\Gamma$,
and $w=w(x,t)\in\R^m$, $v_0\in\R^m$ are suitable vectors.
In $n=3$ dimensions, it holds
\begin{align}
\label{A0def}
A_0(x,t) =\;& \left(\begin{array}{c;{2pt/2pt}c;{2pt/2pt}c;{2pt/2pt}c;{2pt/2pt}c}
\C_1^1 & \C_2^1 & \C_3^1 & 0 & 0\\ \hdashline[2pt/2pt]
\C_1^2 & \C_2^2 & \C_3^2 & 0 & 0\\ \hdashline[2pt/2pt]
\C_1^3 & \C_2^3 & \C_3^3 & 0 & 0\\ \hdashline[2pt/2pt]
0 & 0 & 0 & \rho\,\id_{3\times3} & 0\\ \hdashline[2pt/2pt]
0 & 0 & 0 & 0 & \id_{3\times3}\end{array}\right)\!(x,t)\in
\R^{15\times15},\\
\label{Akdef}
A_k(x,t) =\;& \left(\begin{array}{c;{2pt/2pt}c;{2pt/2pt}c;{2pt/2pt}c;{2pt/2pt}c}
0 & 0 & 0 & -\C_1^k & 0\\ \hdashline[2pt/2pt]
0 & 0 & 0 & -\C_2^k & 0\\ \hdashline[2pt/2pt]
0 & 0 & 0 & -\C_3^k & 0\\ \hdashline[2pt/2pt]
-\C_1^k & -\C_2^k & -\C_3^k & 0 & 0\\ \hdashline[2pt/2pt]
0 & 0 & 0 & 0 & \id_{3\times3}
\end{array}\right)\!(x,t)\in\R^{15\times15},\qquad k=1,2,3,
\end{align}
\begin{equation}
\label{Bdef}
B=\left(\!\!\begin{array}{c;{2pt/2pt}c;{2pt/2pt}c;{2pt/2pt}c;{2pt/2pt}c}
0 & 0 & 0 & 0 & 0\\ \hdashline[2pt/2pt]
0 & 0 & 0 & 0 & 0\\ \hdashline[2pt/2pt]
0 & 0 & 0 & 0 & 0\\ \hdashline[2pt/2pt]
{\scriptsize
\begin{array}{@{\!}l@{\;\;\;}l@{\;\;\;}l@{\!\!}}
D_{11}^1 & D_{12}^1 & D_{13}^1\\
D_{11}^2 & D_{12}^2 & D_{13}^2\\
D_{11}^3 & D_{12}^3 & D_{13}^3
\end{array}
}
&
{\scriptsize
\begin{array}{@{\!}l@{\;\;\;}l@{\;\;\;}l@{\!\!}}
D_{21}^1 & D_{22}^1 & D_{23}^1\\
D_{21}^2 & D_{22}^2 & D_{23}^2\\
D_{21}^3 & D_{22}^3 & D_{23}^3
\end{array}
}
&
{\scriptsize
\begin{array}{@{\!}l@{\;\;\;}l@{\;\;\;}l@{\!\!}}
D_{31}^1 & D_{32}^1 & D_{33}^1\\
D_{31}^2 & D_{32}^2 & D_{33}^2\\
D_{31}^3 & D_{32}^3 & D_{33}^3
\end{array}
}
& 
{\scriptsize
\begin{array}{@{\!}l@{\;\;\;}l@{\;\;\;}l@{\!\!}}
\paj S_{1j1}-\pat\rho & \paj S_{1j2} & \paj S_{1j3}\\
\paj S_{2j1} & \paj S_{2j2}-\pat\rho & \paj S_{2j3}\\
\paj S_{3j1} & \paj S_{3j2} & \paj S_{3j3}-\pat\rho
\end{array}
}
& 0\\ \hdashline[2pt/2pt]
-\id & -\id & -\id & -\id & 0 \end{array}\!\!\right).
\end{equation}
The boundary matrix and the vector of the right hand side are given by
\begin{equation}
\label{Mwdef}
M=\left(\begin{array}{c;{2pt/2pt}c;{2pt/2pt}c}
\C_\nu\in\R^{9\times9} & 0\in\R^{9\times3} & 0\in\R^{9\times3}
\\ \hdashline[2pt/2pt]
0\in\R^{3\times9} & 0\in\R^{3\times3} & 0\in\R^{3\times3}\\ \hdashline[2pt/2pt]
0\in\R^{3\times3} & 0\in\R^{3\times3} & \id\in\R^{3\times3}\end{array}\right),
\qquad w(x,t)=
{\scriptsize
\left(\begin{array}{c}
\vec{0}\in\R^9\\ \hdashline[2pt/2pt]
{\scriptsize \rho(x,t)\,e_1(x,t)}\\ \rho(x,t)\,e_2(x,t)\\ \rho(x,t)\,e_3(x,t)
\\ \hdashline[2pt/2pt]
\vec{0}\in\R^3 \end{array}\right)
}.
\end{equation}
\end{Lemma}

\vspace*{2mm}
\noindent {\bf Proof.}
The following procedure is a modification of an example in
\cite[p.163]{John82} for a scalar hyperbolic equation.
Subsequently we restrict to the case $n=3$ which
allows us to explicitly write down all matrices and explain the necessary
transformations. However, the method is valid for any dimension $n\in\N$.
Let for $n=3$
\begin{align}
v=\,&(v_1,\ldots,v_{15})^T\in\R^{15}\nn\\
\label{vdef}
:=\,&(\partial_1\ut_1,\partial_1\ut_2,\partial_1\ut_3,\partial_2\ut_1,
\partial_2\ut_2,\partial_2\ut_3,\partial_3\ut_1,\partial_3\ut_2,
\partial_3\ut_3,\pat\ut_1,\pat\ut_2,\pat\ut_3,\ut_1,\ut_2,\ut_3)^T.
\end{align}
Note that the last 3 components $\ut_1$, $\ut_2$ and $\ut_3$ of $v$ are only
necessary to incorporate the Dirichlet boundary conditions (\ref{mod3}).
We have a first set of compatibility equations
\begin{align}
\pat v_1-\partial_1 v_{10} =&\; 0,\qquad \pat v_4-\partial_2 v_{10}=0,\qquad
\pat v_7-\partial_3 v_{10}=0,\nn\\
\label{comp1}
\pat v_2-\partial_1 v_{11} =&\; 0,\qquad \pat v_5-\partial_2 v_{11}=0,\qquad
\pat v_8-\partial_3 v_{11}=0,\\
\pat v_3-\partial_1 v_{12} =&\; 0,\qquad \pat v_6-\partial_2 v_{12}=0,\qquad
\pat v_9-\partial_3 v_{12}=0.\nn
\end{align}
The equations (\ref{mod}) determine $\pat v_{10}$, $\pat v_{11}$ and
$\pat v_{12}$. For the last variables, there is a second set of compatibility
relations
\begin{equation}
\label{comp2}
\pat v_{13}-v_{10}=0,\qquad \pat v_{14}-v_{11}=0,\qquad \pat v_{15}-v_{12}=0
\end{equation}
together with
\begin{align}
\partial_1 v_{13}-v_1 =&\; 0,\qquad \partial_2 v_{13}-v_4=0,\qquad
\partial_3 v_{13}-v_7=0,\nn\\
\label{comp3}
\partial_1 v_{14}-v_2 =&\; 0,\qquad \partial_2 v_{14}-v_5=0,\qquad
\partial_3 v_{14}-v_8=0,\\
\partial_1 v_{15}-v_3 =&\; 0,\qquad \partial_2 v_{15}-v_6=0,\qquad
\partial_3 v_{15}-v_9=0.\nn
\end{align}
With (\ref{vdef}), we write the modified linear elasticity equation
\begin{equation}
\label{lelast}
\rho(x,t)\,\patt\ut_i-\C_{ijkl}(x,t)\,\paj\pak\ut_l
+(S_{jki}-S_{ijk})\paj\pat\ut_k=0
\end{equation}
and (\ref{comp1}), (\ref{comp2}), (\ref{comp3}) in the matrix form
\[ \tilde{A}_0\,\pat v+\sum_{k=1}^3\tilde{A}_k\,\pak v+Bv=0. \]
Direct inspection yields
\begin{equation}
\label{Atdef}
\tilde{A}_0(x,t)=\left(\begin{array}{c;{2pt/2pt}c;{2pt/2pt}c}
\id\in\R^{9\times9} & 0\in\R^{9\times3}& 0\in\R^{9\times3}
\\ \hdashline[2pt/2pt]
0\in\R^{3\times9} & \rho(x,t)\,\id\in\R^{3\times3}& 0\in\R^{3\times3}
\\ \hdashline[2pt/2pt]
0\in\R^{3\times9} & 0\in\R^{3\times3} & \id\in\R^{3\times3}
\end{array}\right)\in\R^{15\times15},
\end{equation}
\begin{equation}
\label{At1def}
\tilde{A}_1=\left(\begin{array}{l;{2pt/2pt}l;{2pt/2pt}l}
0\in\R^{3\times9} & -\id\in\R^{3\times3} & 0\in\R^{3\times3}
\\ \hdashline[2pt/2pt]
0\in\R^{3\times9} & 0\in\R^{3\times3} & 0\in\R^{3\times3}
\\ \hdashline[2pt/2pt]
0\in\R^{3\times9} & 0\in\R^{3\times3} & 0\in\R^{3\times3}
\\ \hdashline[2pt/2pt]
{\scriptsize
\begin{array}{@{\!\!}l@{\;}l@{\;}l@{\;}l@{\;}l@{\;}l@{\;}l@{\;}l@{\;}l@{\!\!}}
-\C_{1111}& -\C_{1112}& -\C_{1131}& -\C_{1211}& -\C_{1221}&
-\C_{1231}& -\C_{1311}& -\C_{1321}& -\C_{1331}\\
-\C_{2111}& -\C_{2112}& -\C_{2131}& -\C_{2211}& -\C_{2221}&
-\C_{2231}& -\C_{2311}& -\C_{2321}& -\C_{2331}\\
-\C_{3111}& -\C_{3112}& -\C_{3131}& -\C_{3211}& -\C_{3221}&
-\C_{3231}& -\C_{3311}& -\C_{3321}& -\C_{3331}
\end{array}
}
& F_1\in\R^{3\times3} & 0\in\R^{3\times3}\\ \hdashline[2pt/2pt]
0\in\R^{3\times9} & 0\in\R^{3\times3} & \id\in\R^{3\times3}
\end{array}\right),
\end{equation}
\begin{equation}
\label{At2def}
\tilde{A}_2=\left(\begin{array}{l;{2pt/2pt}l;{2pt/2pt}l}
0\in\R^{3\times9} & 0\in\R^{3\times3} & 0\in\R^{3\times3}
\\ \hdashline[2pt/2pt]
0\in\R^{3\times9} & -\id\in\R^{3\times3} & 0\in\R^{3\times3}
\\ \hdashline[2pt/2pt]
0\in\R^{3\times9} & 0\in\R^{3\times3} & 0\in\R^{3\times3}
\\ \hdashline[2pt/2pt]
{\scriptsize
\begin{array}{@{\!\!}l@{\;}l@{\;}l@{\;}l@{\;}l@{\;}l@{\;}l@{\;}l@{\;}l@{\!\!}}
-\C_{1112}& -\C_{1122}& -\C_{1132}& -\C_{1212}& -\C_{1222}&
-\C_{1232}& -\C_{1312}& -\C_{1322}& -\C_{1332}\\
-\C_{2112}& -\C_{2122}& -\C_{2132}& -\C_{2212}& -\C_{2222}&
-\C_{2232}& -\C_{2312}& -\C_{2322}& -\C_{2332}\\
-\C_{3112}& -\C_{3122}& -\C_{3132}& -\C_{3212}& -\C_{3222}&
-\C_{3232}& -\C_{3312}& -\C_{3322}& -\C_{3332}
\end{array}
}
& F_2\in\R^{3\times3} & 0\in\R^{3\times3}\\ \hdashline[2pt/2pt]
0\in\R^{3\times9} & 0\in\R^{3\times3} &
\id\in\R^{3\times3}\end{array}\right),
\end{equation}
\begin{equation}
\label{At3def}
\tilde{A}_3=\left(\begin{array}{l;{2pt/2pt}l;{2pt/2pt}l}
0\in\R^{3\times9} & 0\in\R^{3\times3} & 0\in\R^{3\times3}
\\ \hdashline[2pt/2pt]
0\in\R^{3\times9} & 0\in\R^{3\times3} & 0\in\R^{3\times3}
\\ \hdashline[2pt/2pt]
0\in\R^{3\times9} & -\id\in\R^{3\times3} & 0\in\R^{3\times3}
\\ \hdashline[2pt/2pt]
{\scriptsize
\begin{array}{@{\!\!}l@{\;}l@{\;}l@{\;}l@{\;}l@{\;}l@{\;}l@{\;}l@{\;}l@{\!\!}}
-\C_{1113}& -\C_{1123}& -\C_{1233}& -\C_{1213}& -\C_{1223}&
-\C_{1233}& -\C_{1313}& -\C_{1323}& -\C_{1333}\\
-\C_{2113}& -\C_{2123}& -\C_{2133}& -\C_{2213}& -\C_{2223}&
-\C_{2233}& -\C_{2313}& -\C_{2323}& -\C_{2333}\\
-\C_{3113}& -\C_{3123}& -\C_{3233}& -\C_{3213}& -\C_{3223}&
-\C_{3233}& -\C_{3313}& -\C_{3323}& -\C_{3333}
\end{array}
}
& F_3\in\R^{3\times3} & 0\in\R^{3\times3}\\ \hdashline[2pt/2pt]
0\in\R^{3\times9} & 0\in\R^{3\times3} &
\id\in\R^{3\times3}\end{array}\right)
\end{equation}
for the $3\times3$-matrices
{\small
\begin{equation}
\label{Fidef}
\begin{aligned}
F_1:=\;& \left(\!\!\begin{smallmatrix}
0 & S_{121}-S_{112}\; & S_{131}-S_{113}\\
S_{112}-S_{211}\; & S_{122}-S_{212}\; & S_{132}-S_{213}\\
S_{113}-S_{311}\; & S_{123}-S_{312}\; & S_{133}-S_{313}
\end{smallmatrix}\!\!\right),\quad
F_2:=\; \left(\!\!\begin{smallmatrix}
S_{211}-S_{121}\; & S_{221}-S_{122}\; & S_{231}-S_{123}\\
S_{212}-S_{221}\; & 0\; & S_{232}-S_{223}\\
S_{213}-S_{321}\; & S_{223}-S_{322}\; & S_{233}-S_{323}
\end{smallmatrix}\!\!\right),\\
F_3:=\;& \left(\!\!\begin{smallmatrix}
S_{311}-S_{131}\; & S_{321}-S_{132}\; & S_{331}-S_{133}\\
S_{312}-S_{231}\; & S_{322}-S_{232}\; & S_{332}-S_{233}\\
S_{313}-S_{331}\; & S_{323}-S_{332}\; & 0
\end{smallmatrix}\!\!\right).
\end{aligned}
\end{equation}
}
However, $\tilde{A}_1(x,t)$, $\tilde{A}_2(x,t)$,
$\tilde{A}_3(x,t)\in\R^{15\times15}$ are not symmetric.
The matrices $F_i$ are on the diagonal of $\tilde{A}_i$, $i=1,2,3$.
They are not symmetric and cannot be rearranged without destroying the
symmetry and positive definiteness of $A_0$. Imposing the strong
assumption~(\ref{Ssym2b}), we obtain $F_1=F_2=F_3=0\in\R^{3\times3}$.

Now, to symmetrize $\tilde{A}_i$, we form suitable linear
combinations of the compatibility equations (\ref{comp1}).
This idea is first presented in \cite{Sfyris24}.

Exemplary, to get the first line of $A_1$, we use (cf. the first column of
$\tilde{A}_1$)
\begin{equation}
\label{lc11}
\C_{1111}(\pat v_1-\partial_1 v_{10})+\C_{2111}(\pat v_2-\partial_1 v_{11})
+\C_{3111}(\pat v_3-\partial_1 v_{12})=0,
\end{equation}
to get the second line of $A_1$, we use (cf. the second column of
$\tilde{A}_1$)
\begin{equation}
\label{lc12}
\C_{1112}(\pat v_1-\partial_1 v_{10})+\C_{2112}(\pat v_2-\partial_1 v_{11})
+\C_{3112}(\pat v_3-\partial_1 v_{12})=0,
\end{equation}
and eventually to get the ninth line of $A_1$, we use
(cf. the ninth column of $\tilde{A}_1$)
\begin{equation}
\label{lc19}
\C_{1331}(\pat v_1-\partial_1 v_{10})+\C_{2331}(\pat v_2-\partial_1 v_{11})
+\C_{3331}(\pat v_3-\partial_1 v_{12})=0.
\end{equation}
So we obtain $A_1$, see (\ref{Akdef}) with $k=1$.
In the same way, $\tilde{A}_2$ is symmetrized. For the first line of $A_2$,
we use (cf. 1. column of $\tilde{A}_2$)
\begin{equation}
\label{lc21}
\C_{1112}(\pat v_4-\partial_2 v_{10})+\C_{2112}(\pat v_5-\partial_2 v_{11})
+\C_{3112}(\pat v_6-\partial_2 v_{12})=0,
\end{equation}
for the second line
\begin{equation}
\label{lc22}
\C_{1122}(\pat v_4-\partial_2 v_{10})+\C_{2122}(\pat v_5-\partial_2 v_{11})
+\C_{3122}(\pat v_6-\partial_2 v_{12})=0,
\end{equation}
and so forth. Finally, to symmetrize $\tilde{A}_3$, we use for the first
two lines
\begin{align}
\label{lc31}
\C_{1113}(\pat v_7-\partial_3 v_{10})+\C_{2113}(\pat v_8-\partial_3 v_{11})
+\C_{3113}(\pat v_9-\partial_3 v_{12})=&\; 0,\\
\label{lc32}
\C_{1123}(\pat v_7-\partial_3 v_{10})+\C_{2123}(\pat v_8-\partial_3 v_{11})
+\C_{3123}(\pat v_9-\partial_3 v_{12})=&\; 0,
\end{align}
and similar operations for lines 3 to 9. With these operations we obtain
$A_2$, $A_3$, see (\ref{Akdef}) with $k=2,3$.
The linear combinations (\ref{lc11})--(\ref{lc19}), (\ref{lc21})--(\ref{lc32})
lead as well to changes in $\tilde{A}_0$, resulting in (\ref{A0def}).
For instance, (\ref{lc11}) modifies the first line of $A_0$ which
results in three non-zero entries, (\ref{lc12}) changes the second line of
$A_0$. In total, the coefficients in (\ref{lc11})--(\ref{lc19}) constitute
the first $3\times3$-block $\C_1^1$ in $A_0$.

\vspace*{2mm}
The major and minor symmetry (\ref{Csym}) of $\C$ imposes the symmetry of
the matrices $\C_j^i$ in (\ref{Cijdef}) as well as $\C_j^i=\C_i^j$ for
$i,j\in\{1,2,3\}$, implying the symmetry of $A_0(x,t)$ for all $(x,t)$.
Due to $\rho(x,t)>0$ in $\cD$ and the positive definiteness of $\C$,
$A_0(x,t)$ is positive definite.  The components $-v_{10}$ in (\ref{comp2})$_1$,
$-v_{11}$ in (\ref{comp2})$_2$ and $-v_{12}$ in (\ref{comp2})$_3$ lead to the
last $-\id\in\R^{3\times3}$ block of $B$ in (\ref{Bdef}).
The other $-\id$ blocks in $B$ are due to the terms $-v_1$, \ldots,
$-v_{9}$ in (\ref{comp3}). The $\id\in\R^{3\times3}$ block in
(\ref{At1def})--(\ref{At3def}) is due to $\partial_1v_{13}$, \ldots,
$\partial_3v_{15}$ in (\ref{comp3}).
This completes the reformulation of (\ref{comp1})--(\ref{lelast})
in form of a linear symmetric first-order hyperbolic system.

The remaining terms of (\ref{mod}) can be incorporated in $B$, leading to
(\ref{Bdef}). The boundary matrix $M$ and the right hand side $w$
are specified by (\ref{Mwdef}).

We want to comment on the matrix $\C_\nu$ appearing in $M$. Due to Dirichlet
boundary conditions, it holds $\ut=0$ on $\Gamma$. As a consequence,
all tangential derivatives of $\ut$ must vanish along $\Gamma$.
Exemplary, consider the special
case of the local coordinates introduced before Definition~\ref{def2} below,
where $\Gamma=\partial\Omega$ is straight and corresponds to $x_1=0$.
Here, the tangential derivatives $\partial_2\ut_i$ and $\partial_3\ut_i$,
$i\in\{1,2,3\}$ vanish on $\Gamma$. In matrix form this reads
\begin{equation}
\label{Pmat}
P(\!\!
\begin{array}{c@{\,}c@{\,}c@{\,\vdots\,}c@{\,}c@{\,}c@{\,\vdots\,}c@{\,}c@{\,}c}
v_1 & v_2 &v_3 & v_4 & v_5 &v_6 & v_7 & v_8 &v_9 \end{array}\!\!)^T
:=\left(\begin{array}{c;{2pt/2pt}c;{2pt/2pt}c}
0 & 0 & 0\\ \hdashline[2pt/2pt]
0 & \id & 0\\ \hdashline[2pt/2pt]
0 & 0 & \id\end{array}\right)(\!\!
\begin{array}{c@{\,}c@{\,}c@{\,\vdots\,}c@{\,}c@{\,}c@{\,\vdots\,}c@{\,}c@{\,}c}
v_1 & v_2 &v_3 & v_4 & v_5 &v_6 & v_7 & v_8 &v_9 \end{array}\!\!)^T=0.
\end{equation}
In symmetrizing $\tilde{A}_k$ for $k=1,2,3$, the linear combinations
(\ref{lc11})--(\ref{lc32}) are applied, converting the matrix
$P\in\R^{9\times9}$ to $\C_\nu$ with $\nu=(-1,0,0)^T$ in the example.
This is in line with the transformation of the upper left
$9\times9$ block $\id$ of $\tilde{A}_0$ to the upper left $9\times9$ block of
$A_0$ in (\ref{A0def}).

\vspace*{2mm}
The vector $v_0$ for the initial values can be directly read off from
(\ref{mod1}), (\ref{mod2}),
{\scriptsize
\begin{align}
\label{v0def}
v_0(x) =\Big(& \partial_1 u_{0,1}(x)-\partial_1\overline{u}_1(x,0),
\partial_1 u_{0,2}(x)-\partial_1\overline{u}_2(x,0),
\partial_1 u_{0,3}(x)-\partial_1\overline{u}_3(x,0),\ldots,
\partial_3 u_{0,3}(x)-\partial_3\overline{u}_3(x,0),\\
& \frac{\mu_{01}(x)-S_{1kl}(x,0)\pak(u_{0l}(x)+\ovu_l(x,0))}{\rho_0(x)}
-\pat\overline{u}_1(x,0), \frac{\mu_{02}(x)-S_{2kl}(x,0)
\pak(u_{0l}(x)+\ovu_l(x,0))}{\rho_0(x)}-\pat\overline{u}_2(x,0),\nn\\
& \frac{\mu_{03}(x)-S_{3kl}(x,0)\pak(u_{0l}(x)+\ovu_l(x,0))}
{\rho_0(x)}-\pat\overline{u}_3(x,0),u_{0,1}(x)-\overline{u}_1(x,0),
u_{0,2}(x)-\overline{u}_2(x,0),u_{0,2}(x)-\overline{u}_2(x,0)\big)^T.\nn
\end{align}
}
With (\ref{A0def}), (\ref{Akdef}), (\ref{Bdef}) and (\ref{Mwdef}),
the equivalence of (\ref{syst1})--(\ref{syst3}) with
(\ref{mod})--(\ref{mod3}) has been shown. \qed

\section{Existence and uniqueness of weak solutions}
\label{secexist}
In this section we apply the existence theory for linear symmetric hyperbolic
systems of first order to the Willis system. A good general
introduction and overview of mathematical methods for hyperbolic systems can be
found in \cite[Chapter~7]{Evans19}. An early $L^2$-theory for linear symmetric
hyperbolic systems in bounded domains is developed in \cite{Friedrichs58},
see also \cite{Friedrichs54}.
The case where the boundary is non-characteristic (see Definition~\ref{def1}
below for explanations) is covered in \cite{RauchMassey} and
\cite{LaxPhillips}.
In the situation studied here, $\Gamma$ is characteristic of constant
multiplicity. This has been analysed for tangential regularity in \cite{Rauch85}
and more generally in \cite{Shizuta95}. We also learned a lot from
the seminal paper \cite{HM77}. An alternative approach would be the use of
semigroup theory, see, e.g. \cite{Xin09}, which however appears to be less
flexible.

\subsection{Existence theory for $\Omega=\R^n$}
As long as no boundary conditions are involved, the proof of existence and
uniqueness of solutions to linear symmetric hyperbolic first-order systems
is straightforward and we begin with this case.
The following theorem is taken from \cite{HM77}.
\begin{Theorem}[Existence and uniqueness for $\Omega=\R^n$]
\label{theo1}
Consider the linear symmetric first-order hyperbolic system
(\ref{linsym}) on $\R^n$ with initial data $v_0$.
Let $s\in\N$ and assume that
\begin{itemize}
\item[\rm(i)] $A_0$, $A_i$ and $B$ are in
$C_b^\infty(\R^n\times[0,T];\,\R^{m\times m})$.
\vspace*{-2mm}
\item[\rm(ii)] $A_0$ and $A_i$, $1\le i\le n$ are symmetric.
\vspace*{-2mm}
\item[\rm(iii)] $A_0$ is {\em uniformly positive definite}, i.e. there exists a
constant $\delta>0$ with
\begin{equation}
\label{Apos}
\langle\xi,A_0(x,t)\xi\rangle \ge\delta\;\|\xi\|^2\qquad\mbox{for all }\xi\in
\R^m\!\setminus\!\{0\}\mbox{ and all }x\in\R^n,\,t\in[0,T].
\end{equation}
\vspace*{-8mm}
\item[\rm(iv)] $w\in H^s(\R^n\times[0,T];\,\R^m)$.
\vspace*{-2mm}
\item[\rm(v)] $v_0\in H^s(\R^n;\,\R^m)$.
\end{itemize}
Then there exists a unique solution $v$ of (\ref{linsym}) in $\R^n$ belonging to
$C^r([0,T];\,H^{s-r}(\R^n;\,\R^m))$ for $0\le r\le s$, such that
$v(\cdot,0)=v_0$.
The solution varies continuously with the initial data in $H^s(\R^n;\,\R^m)$.
Finally, the equations are hyperbolic in the sense that if $v_0$ and $w$ have
compact support then so does $v(\cdot,t)$ for each $t$.
\end{Theorem}

An immediate consequence is
\begin{Corollary}
\label{Cor1}
Let the assumptions (A0)--(A2) and (\ref{Ssym2}) hold, and let
(A3) be satisfied for an integer $s\ge1$. Then there exists a unique solution
\begin{equation}
\label{vreg}
v\in\bigcap_{r=0}^sC^r([0,T];\,H^{s-r}(\R^n;\,R^m))
\end{equation}
to the symmetric linear hyperbolic system (\ref{syst1}), (\ref{syst3}) with
$\Omega=\R^n$.
Consequently, there exists a unique solution vector $(\mu,\sigma,u)$ to
(\ref{S1o})--(\ref{BC2}) of the Willis system in $\Omega=\R^n$ satisfying
\begin{align}
\label{ureg1}
u\in\;& \bigcap_{r=0}^{s+1}C^r([0,T];\,H^{s+1-r}(\Omega;\,\R^n)),\\
\label{mureg}
\mu,\,\sigma\in\;& \bigcap_{r=0}^sC^r([0,T];\,H^{s-r}(\Omega;\,\R^n)).
\end{align}
\end{Corollary}

\noindent{\bf Proof.} First we verify that the matrices $A_0$, $A_i$ and $B$
introduced in Lemma~\ref{lem1}, (\ref{A0def})--(\ref{Bdef})
satisfy the requirements (i)--(iii) of Theorem~\ref{theo1}.
Evidently, (i) follows from (A0), (A1) and (A2).
The symmetry (ii) is a direct consequence of (\ref{Csym}) and (\ref{A0def}),
(\ref{Akdef}).
By the uniform positive definiteness of $\C$, and Sylvester's criterion,
all principal minors of $\C$ are strictly positive.
Together with the positivity condition (\ref{rhopos}) on $\rho$,
this yields that $A_0$ is uniformly positive definite.
In the case $\Omega=\R^n$, one can formally set $\ovu:=0$ in the derivation of
the system leading to $u=\ut$ and $w(x,t)\equiv0$ such that (iv) holds.
Finally, due to (\ref{v0def}) and (A0)--(A3),
we obtain $v_0\in H^s(\R^n;\,\R^m)$. Hence, (i) to (v) of
Theorem~\ref{theo1} are satisfied, yielding the existence of a unique solution
$v\in W:=\bigcap_{r=0}^sC^r([0,T];\,H^{s-r}(\R^n;\,\R^m))$.
By (\ref{vdef}), it holds $\ve(u)$, $Du$, $\pat u\in W$, implying
$u\in\bigcap_{r=0}^sC^r([0,T];\,H^{s+1-r}(\R^n;\,\R^m))$ and
$u\in\bigcap_{r=0}^sC^{r+1}([0,T];\,H^{s-r}(\R^n;\,\R^m))$. This can be
subsumed in one formula (\ref{ureg1}).
With the smoothness (\ref{CCb}), (\ref{SCb}) of $\C$ and $S$, Eq.~(\ref{S2o})
yields $\sigma\in W$. With the smoothness (\ref{rhoCb}) of $\rho$ and
(\ref{SCb}) of $S$, Eq.~(\ref{S3o}) yields $\mu\in W$. This proves
(\ref{mureg}). \qed

\begin{Remark}
\label{rem3}
We use the embedding
\[ W^{s,p}(\R^n)\hookrightarrow C^{r,1-\frac{n}{p}}(\R^n) \]
with $ps>n$, $r+\alpha=s-n/p$ for $\alpha\in(0,1)$,
see \cite[Theorem~9.12]{Brezis11}, setting $p=2$.

\noindent In $n=3$ dimensions we obtain the continuity of the embeddings
\begin{equation}
\label{embed}
H^2(\R^3)\hookrightarrow C^{0,\frac12}(\R^3),\quad
H^3(\R^3)\hookrightarrow C^{1,\frac12}(\R^3),\quad
H^4(\R^3)\hookrightarrow C^{2,\frac12}(\R^3).
\end{equation}
For $\mu_0\in H^2(\R^3;\,\R^3)$, $u_0\in H^3(\R^3;\,\R^3)$, (\ref{mureg}) and
(\ref{embed}) yield
\begin{align}
\label{u3da}
u\in\;& C^0([0,T];\,C^{1,1/2}(\R^3;\,\R^3))\cap
C^1([0,T];\,C^{0,1/2}(\R^3;\,\R^3)),\\
\label{mu3da}
\mu,\,\sigma\in\;& C^0([0,T];\,C^{0,1/2}(\R^3;\,\R^3))\cap
C^1([0,T];\,H^1(\R^3;\,\R^3)),
\end{align}
while for $\mu_0\in H^3(\R^3;\,\R^3)$, $u_0\in H^4(\R^3;\,\R^3)$ we even obtain
\begin{align}
\label{u3db}
& u\!\in\! C^0([0,T];C^{2,1/2}(\R^3;\R^3))\cap
C^1([0,T];C^{1,1/2}(\R^3;\R^3))\cap C^0([0,T];C^{0,1/2}(\R^3;\R^3)),\\
\label{mu3db}
& \mu,\,\sigma\in C^0([0,T];C^{1,1/2}(\R^3;\R^3))\cap
C^1([0,T];C^{0,1/2}(\R^3;\R^3))
\end{align}
for the unique solution vector $(\mu,\sigma,u)$ of (\ref{S1o})--(\ref{BC2}).
\end{Remark}

\vspace*{2mm}
\subsection{Existence theory in a bounded domain}
Now we turn to the boundary value problem (\ref{syst1})--(\ref{syst3}) for
a bounded domain $\Omega\subset\R^n$.

\begin{Corollary}[Classical Dirichlet boundary conditions]
\label{Cor2}
Let $\Omega\subset\R^n$ be a bounded domain and assume
\begin{equation}
\label{ovusupp}
\supp(u_0),\,\supp(\mu_0),\,\supp(\ovu_i(t)),\, \supp(\paj\ovu_i(t))
\subset\Omega\qquad\mbox{for }t\in[0,T],\; i,j\in I.
\end{equation}
Let the assumptions (A0)--(A2) and (\ref{Ssym2}) hold and let (A3) be satisfied
for an integer $s\ge1$. Then, Corollary~\ref{Cor1} remains true, i.e. there
exists a unique solution $(\mu,\sigma,u)$ to (\ref{S1o})--(\ref{BC3}) and
(\ref{vreg})--(\ref{mureg}) hold for a bounded domain $\Omega$.
In $n=3$ dimensions, if\/ $\Omega$ has Lipschitz boundary, for
$\mu_0\in H^2(\Omega;\,\R^3)$, $u_0\in H^3(\Omega;\,\R^3)$, it holds
\begin{align}
\label{uom}
u\in\;& C^0([0,T];\,C^{1,1/2}(\Omega;\,\R^3))\cap
C^1([0,T];\,C^{0,1/2}(\Omega;\,\R^3)),\\
\label{muom}
\mu,\,\sigma\in\;& C^0([0,T];\,C^{0,1/2}(\Omega;\,\R^3)),\quad
\mu\in C^1([0,T];\,H^1(\Omega;\,\R^3))
\end{align}
for the unique solution vector $(\mu,\sigma,u)$ of (\ref{S1o})--(\ref{BC3}).
\end{Corollary}

\noindent{\bf Proof.} The condition~(\ref{ovusupp}) implies
$\supp(\pat\ovu)$, $\supp(\patt\ovu)\subset\Omega$ and 
$\supp(\paj u_{0})\subset\Omega$.
So, by (\ref{edef}), $\supp((\rho e_i)(t))\subset\Omega$ for $t\in[0,T]$ and
$i\in\{1,2,3\}$. Using (\ref{Mwdef}), (\ref{v0def}) and (\ref{ovusupp}),
we find
\[ \supp(v_0),\,\supp(w(t))\subset\Omega\qquad\mbox{for }t\in[0,T]. \]
With Theorem~\ref{theo1}, this ensures $\mathrm{supp}\,(v(t))\subset\Omega$ for
$t\in[0,T]$. If $\partial\Omega$ is Lipschitz continuous, there is a compact
embedding $H^s(\Omega)\hookrightarrow C^{0,1/2}(\overline{\Omega})$ for
$s>3/2$, see \cite[Theorem~6 p.270]{Evans19},
leading to (\ref{uom}), (\ref{muom}). \qed

\vspace*{2mm}
Note that (\ref{ovusupp}) implies $\ovu=0$ on $\partial\Omega$, i.e.
classical Dirichlet boundary conditions.

\vspace*{5mm}
Now we discuss general, possibly non-regular Dirichlet boundary conditions.
The following analysis is based on the methods developed in
\cite{Shizuta95} and \cite{Rauch85}.

\begin{Definition}
\label{def1}
Let $\nu(x)=(\nu_1,\ldots,\nu_n)(x)$ be the unit outer normal to $\Omega$ in
$x\in\Gamma=\partial\Omega$. The {\em boundary matrix} is given by
\begin{equation}
\label{Anudef}
A_{\nu(x)}:=\sum_{k=1}^n\nu_k(x)A_k(x,t).
\end{equation}
If $A_\nu$ is invertible everywhere on $\Gamma$, the boundary is called
{\em non-characteristic}. If $A_\nu$ is not invertible but has constant rank
on $\Gamma$, the boundary is called
{\em characteristic of constant multiplicity}.
\end{Definition}

A particular difficulty for the theory of boundary value problems is the
possible loss of derivatives in normal direction. For that reason we need to
introduce new spaces.

Let $(\mathcal{O}_i,\chi_i)_{1\le i\le l}$ be a partition of unity of $\Gamma$.
For some $\delta>0$, let
\[ \Omega_\delta:=\{x\in\Omega\;|\;\mathrm{dist}(x,\Gamma)>\delta\} \]
and let $\chi_0$ be a smooth function with $\chi_0\!=\!1$ in
$\Omega_\delta$ and $\chi_0\!=\!0$ in a neighborhood of $\Gamma$.\linebreak
Assume $\sum_{i=0}^l\chi_i^2=1$ in $\overline{\Omega}$. Introduce local
coordinates and consider a family of diffeomorphisms
$(\tau_i)_{1\le i\le l}$ from $\R^n$ to $\R^n$ such that $\Gamma$ corresponds
to $x_1=0$ and $\Omega$ corresponds to
\[ B_1^+:=\{x\in\R^n\;|\;|x|<1,\,x_1>0\}. \]
We recall that $\Lambda\in C^\infty(\overline{\Omega};\,\R^n)$ is a
{\it tangential vector field} if
\begin{equation}
\label{tang}
\langle\Lambda(x),\nu(x)\rangle=0\qquad\mbox{for all }x\in\Gamma.
\end{equation}

The following three definitions are taken from \cite{Shizuta95}.
\begin{Definition}
\label{def2}
Let $s\ge0$ be an integer. We introduce $H_*^s(\Omega)$ as
the set of functions $f\in L^2(\Omega;\,\R)$ with the following property:
Let $\Lambda_1$, $\Lambda_2$, \ldots, $\Lambda_j$ be tangential vector fields
and $\Lambda_1'$, $\Lambda_2'$, \ldots, $\Lambda_k'$ be non-tangential vector
fields. Then
\begin{equation}
\label{Hsdef}
\Lambda_1\Lambda_2\ldots\Lambda_j\Lambda_1'\Lambda_2'\ldots\Lambda_k'f
\in L^2(\Omega)\qquad\mbox{for }j+2k\le s.
\end{equation}
The norm in $H_*^s(\Omega)$ is given by
\begin{equation}
\label{Hsnorm}
\|f\|_{H_*^s(\Omega)}^2:=\|\chi_0f\|_{H^s(\Omega)}^2+\sum_{i=1}^m
\sum_{|\alpha|+2k\le s}\|\partial_\mathrm{tan}^\alpha
\partial_1^kf^{(i)}\|_{L^2(B_1^+)}^2,
\end{equation}
where $f^{(i)}:=(\chi_if)\circ\tau_i^{-1}$, $\alpha=(\alpha_1,\ldots,\alpha_n)$
is a multi-index, $|\alpha|:=\alpha_1+\ldots+\alpha_n$, and
\[ \partial_\mathrm{tan}^\alpha:=(x_1\partial_1)^{\alpha_1}\partial_2^{\alpha_2}
\cdots\partial_n^{\alpha_n}. \]
\end{Definition}
The key to understanding this definition is the observation
that in local coordinates for any point in a neighborhood of $\Gamma$,
$x_1\partial_1$, $\partial_2$, \ldots, $\partial_n$ span the tangential vector
fields. The normal vector field $\partial_\nu$ corresponds to $-\partial_1$ in
local coordinates. The space $H_*^1(\Omega)$ is identical to
$H^1_\mathrm{tan}(\Omega)$ introduced in \cite{Rauch85}.
\begin{Definition}
\label{def3}
Let $s\ge0$ be an integer. We introduce $H_{**}^s(\Omega)$
as the set of all functions $f\in L^2(\Omega;\,\R)$ with the following property:
Let $\Lambda_1$, $\Lambda_2$, \ldots, $\Lambda_j$ be tangential vector fields
and let $\Lambda_1'$, $\Lambda_2'$, \ldots, $\Lambda_k'$ be non-tangential
vector fields. Then
\begin{equation}
\label{Hssdef}
\Lambda_1\Lambda_2\cdots\Lambda_j\Lambda_1'\Lambda_2'\cdots\Lambda_k'f
\in L^2(\Omega)\qquad\mbox{for }j+2k\le s+1\mbox{ and }j+k\le m.
\end{equation}
The norm on $H_{**}^s(\Omega)$ is given by
\begin{equation}
\label{Hssnorm}
\|f\|_{H_{**}^s(\Omega)}^2:=\|\chi_0f\|_{H^s(\Omega)}^2
+\sum_{i=1}^l\sum_{\substack{|\alpha|+2k\le s+1\\|\alpha|+k\le s}}
\|\partial_\mathrm{tan}^\alpha\partial_1^kf^{(i)}\|_{L^2(B_1^+)}^2.
\end{equation}
\end{Definition}
In (\ref{Hsnorm}), (\ref{Hssnorm}), different choices of
$(\mathcal{O}_i,\chi_i)_{0\le i\le l}$ and $(\tau_i)_{1\le i\le l}$ lead to
equivalent norms. For any $s\ge0$, there is a continuous embedding
$H^s(\Omega)\hookrightarrow H_{**}^s(\Omega)\hookrightarrow H_*^s(\Omega)$.
Therefore, $H_*^s(\Omega)$ and $H_{**}^s(\Omega)$ may be regarded as subspaces
of $H^s(\Omega)$.

\vspace*{4mm}
\begin{Definition}
\label{def4}
For an integer $s\ge0$, let $X^s([0,T];\,\Omega)$ be the space of functions
$f$ with
\[ \pat^rf\in C^0([0,T];\,H^{m-r}(\Omega)),\qquad 0\le r\le s. \]
Here, $\pat^rf$, $0\le r\le s$ are the derivatives of $f$ in the sense of
distributions. The space $X^s([0,T];\,\Omega)$ is a Banach space with the norm
\begin{align}
\label{Xnorm}
\vertt{f}_{X^s([0,T];\,\Omega)}:=\;& \max_{t\in[0,T]}\vertt{f(t)}_s,\\
\label{xtnorm}
\vertt{f(t)}_s^2 :=\;& \sum_{r=0}^s\|\pat^rf(t)\|_{H^{s-r}(\Omega)}^2.
\end{align}
Analogously, let $X_*^s([0,T];\,\Omega)$ be the space of functions $f$ with
\[ \pat^rf\in C^0([0,T];\,H_*^{s-r}(\Omega)),\qquad 0\le r\le s. \]
The space $X_*^s([0,T];\,\Omega)$ is a Banach space with the norm
\begin{align}
\label{Xsnorm}
\vertt{f}_{X_*^s([0,T];\,\Omega}^2 :=\;& \max_{t\in[0,T]}\vertt{f(t)}_{s,*},\\
\label{Xsnorm1}
\vertt{f(t)}_{s,*}^2 :=\;& \sum_{r=0}^s\|\pat^rf(t)\|_{H_*^{s-r}(\Omega)}^2.
\end{align}
\vspace*{-2mm}
By $W_*^s(0,T;\,\Omega)$ we denote the space of functions $f$ such that
\[ \pat^rf\in L^2(0,T;\,H_*^{s-r}(\Omega)),\qquad 0\le r\le s. \]
For vector-valued functions $f=(f_1,\ldots,f_m)$ with
$f_i\in X^s([0,T];\,\Omega)$ for $1\le i\le m$ we write
$f\in X^s([0,T];\,\Omega)^m$ and $W_*^s(0,T;\,\Omega)^m$ is defined analogously.
\end{Definition}

\vspace*{2mm}
The following theorem is a reformulation of \cite[Theorem~2.1]{Shizuta95}
which also investigates certain quasi-linear problems where $A_0$, $A_1$,
\ldots, $A_n$ may depend on a further function. Theorem~2.1 in
\cite{Shizuta95} is an improved version of Theorem~10 in \cite{Rauch85} which
requires slightly stronger assumptions, only proves tangential regularity,
and does not prove uniqueness of the solution. Note that several
assumptions of Theorem~10 in \cite{Rauch85} are not formulated explicitly,
but are mentioned elsewhere in the article.

\begin{Theorem}
\label{theo2}
Let $s\ge1$ be an integer.
Then an initial boundary value problem (\ref{syst1})--(\ref{syst3}) has a
unique solution $v\in X_*^s([0,T];\,\Omega)^m$ provided the following conditions
are satisfied:
\begin{enumerate}
\item[{\rm(i)}] $\Omega\subset\R^n$ is a bounded open set with boundary $\Gamma$
of class $C^\infty$.
\item[{\rm(ii)}] $M(x)$ is a real matrix-valued function with
$M\in\C^\infty(\Gamma;\,\R^{m\times m})$.
\item[{\rm(iii)}] $A_k(x,t)\in\R^{m\times m}$, $0\le k\le n$ are real symmetric
matrices for every $(x,t)\in\overline{\Omega}\times[0,T]$. The matrix
$A_0(x,t)$ is positive definite for $(x,t)\in\overline{\Omega}\times[0,T]$.
\item[{\rm(iv)}] The dimension of $\cN(x):=\ker A_{\nu(x)}$ and the dimension
of $\ker M(x)$ are constant on each component of $\Gamma$ and it holds
$0<\dim\cN(x)<m$.
\item[{\rm(v)}] $\ker M(x)$ is a maximal nonnegative subspace of
$A_{\nu(x)}$ for $x\in\Gamma$.
\item[{\rm(vi)}] It holds $w\in W_*^s(0,T;\,\Omega)^m$,
$\pat^r w(0)\in H^{s-1-r}(\Omega;\,\R^m)$ for\/ $0\le r\le s-1$ and
the initial values fulfill $v_0\in H^s(\Omega;\,\R^m)$.
\item[{\rm(vii)}] The data $w$, $v_0$ fulfill the compatibility conditions of
order $s-1$ for the initial boundary value problem (\ref{syst1})--(\ref{syst3}).
\end{enumerate}
Then the solution $v$ obeys for positive constants $c_1$, $c_2$ the estimate
\begin{equation}
\label{vreg2}
\vertt{v(t)}_{s,*}\le c_1\Big(\|v_0\|_{H^s(\Omega;\,\R^m)}
+\vertt{w(0)}_{s-1}\Big)\ee^{c_2t}+c_2\int_0^t\ee^{c_2(t-\tau)}
\vertt{w(\tau)}_{s,*}\dtau,\quad t\in[0,T].
\end{equation}
\end{Theorem}

\begin{Remark}
\label{rem4}
The condition (v) in Theorem~\ref{theo2} means that
\begin{equation}
\label{noneg}
\langle A_\nu(x,t)z,z\rangle\ge0\qquad\mbox{for all }(x,t)\in\Gamma\times[0,T],
\,z\in \ker M(x)
\end{equation}
and $\ker M$ as a subspace of $\C^m$
cannot be enlarged while (\ref{noneg}) remains valid.

\vspace*{3mm}
\noindent The compatibility conditions on $w$ and $v_0$ up to order $s-1$ are
canonical as in (\ref{compa}). One takes spatial derivatives up to order $s-1$
of the system, resolves $\pat^pv$ for $0\le p\le s-1$ and evaluates the
result at $t=0$.
In detail, the procedure is as follows, cf. \cite[p.169]{Shizuta95}.
For $p=0$, set $v_{0,0}:=v_0$. For $p=1,2,\ldots, s-1$, set iteratively
\[ v_{0,p}:=\sum_{i=0}^{p-1}
\left(\begin{smallmatrix}p-1\\ i\end{smallmatrix}\right)
G_i(0)v_{0,p-1-i}+\pat^{p\,-\,1}(A_0^{-1}w)(0)\qquad\mbox{in }\Omega, \]
where
\vspace*{-2mm}
\begin{align*}
G_0(t):=\;& -\sum_{k=1}^nA_0^{-1}A_k\pak-A_0^{-1}B,\\
G_i(t):=\;& -\sum_{k=1}^n\pat^i(A_0^{-1}A_k)\pak-\pat^i(A_0^{-1}B),\qquad i\ge1.
\end{align*}
Then, the {\em compatibility conditions} up to order $s-1$ are
\begin{equation}
\label{comps}
Mv_{0,p}=0\qquad\mbox{on }\Gamma\mbox{ for }0\le p\le s-1.
\end{equation}
\end{Remark}
\begin{Corollary}
\label{Cor3}
Let $n=3$ and $\Omega\subset\R^3$ be a bounded open set with boundary $\Gamma$
of class $C^\infty$. Assume (\ref{Ssym2}) and let (A0)--(A4) be
satisfied for an integer $s\ge1$.
Let $u_0$, $\ovu$ satisfy the compatibility conditions up to order $s+1$.
Let $\C_\nu$ be defined by (\ref{Cnudef}) and assume that
\begin{equation}
\label{dimN}
\dim\ker(\C_\nu)=\mathrm{const}>0\qquad\mbox{on every component of }\,\Gamma.
\end{equation}
Then there exists a unique solution $(\mu,\sigma,u)$ of the Willis system
(\ref{S1o})--(\ref{BC2}) which satisfies
\begin{align}
\label{ureg2}
u\in\;& \bigcap_{r=0}^{s+1}C^r([0,T];\,H_*^{s+1-r}(\Omega;\,\R^n)),\\
\label{mureg2}
\mu,\,\sigma\in\;& \bigcap_{r=0}^sC^r([0,T];\,H_*^{s-r}(\Omega;\,\R^n)).
\end{align}
\end{Corollary}

\noindent{\bf Proof.} 
We need to validate (i)--(vii) of Theorem~\ref{theo2} for $n=3$.
Clearly (i) holds by assumption and (ii) follows from (A0) and (\ref{Mwdef}).
Condition~(iii) follows as in Corollary~\ref{Cor1} from Lemma~\ref{lem1}, the
symmetry (\ref{Csym}) and the positivity condition (\ref{rhopos}) on $\rho$.

The condition $v_0\in H^s(\Omega;\,\R^m)$ in (vi) follows from (A0)--(A3).
By (\ref{A4C3a}) and (\ref{edef}) the condition on $\pat^rw(0)$ in
(vi) is ensured. Finally, due to (\ref{A4C3b}), the remaining condition
$w\in W_*^s(0,T;\,\Omega)^m$ from (vi) is also satisfied.

Writing $\C_k^\nu:=\nu_1\C_k^1+\nu_2\C_k^2+\nu_3\C_k^3$ for $k=1,2,3$,
due to (\ref{Akdef}) and (\ref{Anudef}), we have
\begin{equation}
\label{Anu3d}
A_\nu=\left(\!\!\begin{array}{c;{2pt/2pt}c;{2pt/2pt}c;{2pt/2pt}c;{2pt/2pt}c}
0 & 0 & 0 & -\C_1^\nu & 0\\ \hdashline[2pt/2pt]
0 & 0 & 0 & -\C_2^\nu & 0\\ \hdashline[2pt/2pt]
0 & 0 & 0 & -\C_3^\nu & 0\\ \hdashline[2pt/2pt]
-\C_1^\nu & -\C_2^\nu & -\C_3^\nu & 0 & 0\\ \hdashline[2pt/2pt]
0 & 0 & 0 & 0 & (\nu_1+\nu_2+\nu_3)\id
\end{array}\!\!\right)\in\R^{15\times15}.
\end{equation}
Let $z:=(z^1,z^2,\ldots,z^5)^T\in\R^{15}$ be a vector consisting of blocks
$z^k\in\R^3$ for $k=1,\ldots,5$. Evaluating $A_\nu z=0$ leads to
$z^4=z^5=0$ and
\[ \C_1^\nu z^1+\C_2^\nu z^2+\C_3^\nu z^3=0\in\R^3 \]
which is equivalent to $z^4=z^5=0$ and
\begin{equation}
\label{kerCnu}
(z^1,z^2,z^3)^T\in\ker(\C_\nu)
\end{equation}
with $\C_\nu\in\R^{9\times9}$ defined in (\ref{Cnudef}).
So we obtain
\begin{equation}
\label{kerAnu}
\ker A_{\nu(x)}=\Big\{z=(z^1,\ldots,z^5)^T\in\R^{15}\;\Big|\;
z^4=z^5=0,\; (z^1,z^2,z^3)^T\in\ker(\C_{\nu(x)})\Big\},\quad
x\in\Gamma.
\end{equation}
This demonstrates that $0<\dim\ker A_{\nu(x)}<m$ for $x\in\Gamma$.
Due to (\ref{Mwdef}), we find
\begin{equation}
\label{kerM}
\ker M(x)=\Big\{z=(z^1,\ldots,z^5)^T\in\R^{15}\;\Big|\;z^5=0,\;
(z^1,z^2,z^3)^T\in\ker(\C_{\nu(x)})\Big\},\quad x\in\Gamma.
\end{equation}
With Assumption (\ref{dimN}), this implies (iv) of Theorem~\ref{theo2}.

For a vector $z=(z^1,\ldots,z^5)^T\in\ker(M(x))$, the condition
(\ref{noneg}) becomes with $\nu=\nu(x)$
\begin{align*}
\langle A_\nu z,z\rangle =\;& -\langle\C_1^\nu z^4,z^1\rangle-
\langle\C_2^\nu z^4,z^2\rangle-\langle\C_3^\nu z^4,z^3\rangle\\
& -\langle\C_1^\nu z^1,z^4\rangle-\langle\C_2^\nu z^2,z^4\rangle
-\langle\C_3^\nu z^3,z^4\rangle \;\ge\; 0.
\end{align*}
Due to the symmetry of $\C_k^\nu$, this is equivalent to
$\langle\C_1^\nu z^1+\C_2^\nu z^2+\C_3^\nu z^3,z^4\rangle\le0$.
Since $z^4$ for $z=(z^1,z^2,z^3,z^4,z^5)\in\ker(M)$ is arbitrary, we must have
$\C_1^\nu z^1+\C_2^\nu z^2+\C_3^\nu z^3=0$ or
\begin{equation}
\label{noneg2}
\nu_1(\C_1^1z^1+\C_2^1z^2+\C_3^1z^3)+
\nu_2(\C_1^2z^1+\C_2^2z^2+\C_3^2z^3)+
\nu_3(\C_1^3z^1+\C_2^3z^2+\C_3^3z^3)=0.
\end{equation}
This last condition (\ref{noneg2}) is equivalent to $\C_\nu(z^1,z^2,z^3)=0$
which holds due to (\ref{kerM}). Geometrically, Eq.~(\ref{noneg2}) ensures that
all tangential derivatives of $\ut$ vanish along $\Gamma$.
In summary, (\ref{noneg}) holds with equality and enlarging $\ker(M)$
violates (\ref{noneg2}).
This proves the remaining condition (v) in Theorem~\ref{theo2}. \qed

\section{Concluding remarks}
\label{secdis}
In this article, three existence and uniqueness results for (weak) solutions to
a system of partial differential equations related to the Willis model have been
derived, both for the whole space case and for bounded domains. The
investigated system (\ref{S1o})--(\ref{S3o}) differs from the original Willis
equations in that no explicit form of $S$, e.g. no convolution
expression, is postulated.
In addition to the natural symmetry condition (\ref{Ssym}) on $S$ which
guarantees the symmetry of the Cauchy stress tensor $\sigma$, a second
symmetry condition (\ref{Ssym2}) has to be imposed for the analysis.
This condition appears to be necessary mathematically and admits to write
(\ref{S1o})--(\ref{S3o}) as a linear symmetric hyperbolic system of first
order. Combined, (\ref{Ssym}) and (\ref{Ssym2}) impose strong
restrictions on $S$: the third-order Willis coupling tensor must be totally
symmetric. It has to be checked experimentally whether these
conditions are satisfied for a real-world material.
At this point we have no rigorous physical justification for
(\ref{Ssym2}). We refer to \cite{Muhl16} for a discussion of necessary
conditions on the quantities of the Willis system in order to have a physically
correct model. Since inhomogeneous dynamic linear elasticity constitutes a
linear symmetric hyperbolic system, \cite{HM77,Sfyris24}, so should be any
homogenized problem based on the former, here the Willis equations.
The condition~(\ref{Ssym2}) seems to be a welcome
novel restriction on the Willis coupling tensor $S$ introduced by this
requirement.

\vspace*{2mm}
In summary, with the additional symmetry condition (\ref{Ssym2}),
after establishing suitable structural assumptions, mostly
on $S$ and $\rho$, the existence and uniqueness of a (weak) solution follows
from established existence results for linear symmetric first order
hyperbolic systems. If the initial and boundary data is regular
enough, even a unique classical solution is obtained.
Finally, in bounded domains, the condition
(\ref{dimN}) is essential for the existence theory.
Its validity depends crucially on material properties.

\end{document}